\documentclass[a4paper,11pt]{article}
\usepackage{amssymb,amsmath,amsfonts,epsf,amsmath,tikz}
\usepackage{enumerate}
\usepackage{algorithm}
\usepackage{framed}
\usepackage{hyperref}

\tikzset{
  use path/.code={\pgfsyssoftpath@setcurrentpath{#1}}
}
\usetikzlibrary{decorations.pathmorphing}

\newtheorem{theorem}{\bf Theorem}[section]
\newtheorem{corollary}[theorem]{\bf Corollary}
\newtheorem{lemma}[theorem]{\bf Lemma}
\newtheorem{proposition}[theorem]{\bf Proposition}

\newcommand{\proof}{\noindent{\bf Proof.\ }}
\newcommand{\qed}{\hfill $\square$ \medskip}

\newcommand{\Nb}{\mathbb{N}}
\newcommand{\GDG}{D_\mathcal{G}}
\newcommand{\GDR}{D_\mathcal{R}}

\begin{document}

 \title{ Game Distinguishing Numbers of Cartesian Products of Graphs}

\author{
Sylvain Gravier $^{a,b}$ 
\and
Kahina Meslem $^{b,c}$
\and 
Simon Schmidt $^{a,b}$
\and
Souad Slimani $^{a,b,c}$
}

\date{}

\maketitle

\begin{center}
$^a$ Institut Fourier, UMR 5582, Universit\'e Grenoble Alpes\\
100, rue des maths BP74\\ 
38402 Saint-Martin-d'H\`eres Cedex, France\\
\medskip

\medskip
$^b$ S.F.R.~Maths \`a Modeler\\
100, rue des maths BP74\\ 
38402 Saint-Martin-d'H\`eres Cedex, France\\

\medskip
$^c$ LaROMaD, Faculty of Mathematics, U.S.T.H.B.\\
El Alia BP 32 Bab Ezzouar\\
16111 Alger, Algeria\\

\medskip
\url{sylvain.gravier@ujf-grenoble.fr}\\
\url{simon.schmidt@ujf-grenoble.fr}\\
\url{sslimani@usthb.dz}\\
\url{kmeslem@usthb.dz}\\
\medskip

\end{center}
\begin{abstract}
 The distinguishing number of a graph $H$ is a symmetry related graph invariant whose study started two decades ago.
The distinguishing number $D(H)$ is the least integer $d$ such that $H$ has a $d$-distinguishing coloring. A $d$-distinguishing coloring is a
coloring $c:V(H)\rightarrow\{1,\dots,d\}$ invariant only under the trivial automorphism. In this paper, we continue the study of a game variant of this parameter, recently introduced. The distinguishing game is a game with two players, Gentle and Rascal, with antagonist goals. This game is played on a graph $H$ with a fixed set of $d\in\Nb^*$
colors. Alternately, the two players choose a vertex of $H$ and color it with one of the $d$ colors. The game ends when all the vertices have been colored.
Then Gentle wins if the coloring is $d$-distinguishing and Rascal wins otherwise.
This game defines two new invariants, which are the minimum numbers of colors needed to ensure that Gentle has a winning strategy, depending who starts the game.
The invariant could eventually be infinite. In this paper,  we focus on cartesian product, a graph operation
well studied in the classical case. We give sufficient conditions on the order of two connected factors $H$ and $F$ relatively prime, which ensure that one of the game distinguishing numbers of the cartesian product
$H\square F$ is finite. If $H$ is a so-called involutive graph, we give an upper bound of order $D^2(H)$ for one the game distinguishing numbers of $H\square F$. 
Finally, using in part the previous result, we compute the exact value of these invariants for cartesian products of relatively prime cycles. It turns out that the value is either infinite or equal to $2$, depending on the parity of the product order. 
\end{abstract}




\section{Introduction}
In this paper, we consider only simple graphs. For a graph $H$, $V(H)$ and $E(H)$ respectively denote the vertex set and the edge set of $H$.  For an integer $n\geq 3$, $C_n$ is the cycle of order $n$ and for $n\geq 2$, $K_n$ and $P_n$ are respectively the clique and the path of order $n$. The distinguishing number $D(H)$ of a graph $H$ is a symmetry related graph invariant whose study starts two decades ago \cite{albertson}.
More precisely, $D(H)$ is the least integer $d$ such that $H$ has a $d$-distinguishing coloring. A $d$-distinguishing coloring is a vertex-coloring $c:V(H)\rightarrow\{1,\dots,d\}$ invariant only under the trivial automorphism. More generally, we say that an automorphism $\sigma$ of $H$ preserves the coloring $c$ or is a colors
preserving automorphism, if for all $u\in V(H)$, $c(u)=c(\sigma(u))$.
The automorphisms group of $H$ will be denoted by $Aut(H)$.
Clearly, for each coloring $c$ of the vertex set of $H$, the set $Aut_c(H)=\{\sigma \in Aut(H): c\circ\sigma= c\}$ is a subgroup of $Aut(H)$. A coloring $c$ is distinguishing if
$Aut_c(H)$ is trivial.
The ten last years have seen a flourishing number of works on this subject and cartesian products of graphs were thoroughly investigated in \cite{bogstad,klav_power,Imrich_cartes_power,klav.cart.clik,fisher}. In particular, the exact value of $D(K_n\square K_m$) is given in \cite{Imrich_cartes_power,fisher}. 
Another interesting result for our purpose is that if $k\geq 2$, then $D(C_{n_1}\square\cdots\square C_{n_k})=2$, save for $C_3\square C_3$. In that case $D(C_3\square C_3)=3$. This result is an easy consequence of more general results in \cite{Imrich_cartes_power}.
Recently a game variant of the distinguishing number has been introduced in \cite{schmidt}. Defining game invariants for graphs is not a new idea. The two most known game invariants are the game chromatic number, introduced by Brahms in 1981 \cite{faigle}, and the game domination numbers introduced more recently by Bre\v{s}ar, Kla\v{z}ar and Rall \cite{domgame}.

The distinguishing game is a game with two players, Gentle and Rascal, with antagonist goals. This game is played on a graph $H$ with a fixed set of $d\in\Nb^*$
colors. Alternately, the two players choose a vertex of $H$ and color it with one of the $d$ colors. The game ends when all the vertices have been colored. If the coloring is $d$-distinguishing then Gentle wins.  Otherwise Rascal wins.

This game defines two invariants for a graph $H$. The {\em $G$-game distinguishing number} $\GDG(H)$ is the minimum of colors needed to ensure that Gentle has a winning strategy for the game on $H$, assuming he is playing first.
If Rascal is sure to win whatever the number of colors we allow, then $\GDG(H)=\infty$. Similarly, the {\em $R$-game distinguishing number} $\GDR(H)$ is the minimum of colors needed to ensure that Gentle has a winning strategy, assuming Rascal is playing first.
Characterizing graphs with infinite game distinguishing number seems to be a challenging open question. In \cite{schmidt}, the authors give sufficient conditions to have one infinite game distinguishing number.
\begin{proposition} \label{prop:ordertwo}\cite{schmidt}
Let $H$ be a graph and $\sigma$ a non trivial automorphism of $H$ such that $\sigma \circ \sigma=id_{H}$.
\begin{enumerate}
   \item If $|V(H)|$ is even, then $\GDG(H)=\infty$.
   \item If $|V(H)|$ is odd, then $\GDR(H)=\infty$.
\end{enumerate}
\end{proposition}
Also in \cite{schmidt}, the exact values of those invariants have been computed for almost all cycles and hypercubes. And for a large class of graphs, the so-called involutive graphs, a quadratic upper bound involving the classical distinguishing number has been provided (see Section~\ref{sect:inv} for a definition). We give here the precise statement of the results used in this paper. 

\begin{theorem}\label{theo:cycle}\cite{schmidt}
Let $C_n$ be a cycle of order $n\geq 3$.
\begin{enumerate}
\item If $n$ is even (resp. odd), then $\GDG(C_n)=\infty$ (resp. $\GDR(C_n)=\infty$).
 \item If $n$ is even and $n\geq 8$, then $\GDR(C_n)=2$. 
 \item If $n$ is odd, not prime and $n\geq 9$, then $\GDG(C_n)=2$.
 \item If $n$ is prime and $n\geq 5$, then $\GDG(C_n)\leq 3$.
\end{enumerate}
Moreover $\GDR(C_4)=\GDR(C_6)=3$, $\GDG(C_3)=\infty$ and $\GDG(C_5)=\GDG(C_7)=3$.
\end{theorem}
\begin{theorem}\cite{schmidt}\label{theo:finite-inv}
If $H$ is an involutive graph with $D(H)\geq2$, then $\GDR(H)\leq D^2(H)+D(H)-2$.
\end{theorem}  
%
%

In this paper, we deal with cartesian products of connected graphs relatively prime. In Section~\ref{sec:complete}, we prove the following theorem 
which gives sufficient conditions on the order of the two factors to have a finite distinguishing number. 

\begin{theorem}\label{theo:main-cart}
 Let $H$ and $F$ be two non trivial connected graphs relatively prime, with respective order $n$ and $m$. 
\begin{enumerate}
 \item If $n$ is even and $m\geq n-1$, then $\GDR(H\square F)\leq m+1$.
 \item If $n$ is odd, $m$ is even and $m\geq 2n-2$, then $\GDR(H\square F)\leq m+1$.
 \item If $n$ and $m$ are odd and $m\geq 2n-1$, then $\GDG(H\square F)\leq m+1$.
\end{enumerate}
\end{theorem}

In Section~\ref{sect:inv}, we investigate the case where one factor is an involutive graph. In that case, if the classical distinguishing number of the other 
factor is not too big, we have a quadratic upper bound involving the classical distinguishing number of the involutive factor. 

\begin{theorem}\label{theo:main-inv}
 Let $H$ be a connected involutive graph of order $n$, with $D(H)\geq 2$ and $F$ a connected graph relatively prime to $H$. 
If $\displaystyle D(F)\leq \left(\begin{array}{c}
                \frac{n+d^2+d}2-1\\
                \frac{d^2+d}2-1
              \end{array}\right)$, then $\GDR(H\square F)\leq d^2+d-2$, where $d=D(H)$.
\end{theorem}
Finally, in Section~\ref{sec:tore}, we compute the exact value of the two invariants for cartesian products of relatively prime cycles. Since even cycles are involutive graphs, a part of this result 
arises as a corollary of the above theorem.
\begin{theorem}\label{theo:main-tore} 
 Let $n_1,...,n_k$, with $k\geq2$, be $k$ distinct natural numbers greater or equal to $3$.
\begin{enumerate}
 \item If $\displaystyle\prod_{i=1}^k n_i$ is even, then $\GDG(C_{n_1}\square\cdots\square C_{n_k})=\infty$ and $\GDR(C_{n_1}\square\cdots\square C_{n_k})=2$.
 \item If $\displaystyle\prod_{i=1}^k n_i$ is odd, then $\GDR(C_{n_1}\square\cdots\square C_{n_k})=\infty$ and $\GDG(C_{n_1}\square\cdots\square C_{n_k})=2$.
 \end{enumerate}
\end{theorem}
All these three results highly involve the so-called fiber-strategy for Gentle. Section \ref{sect:strat} is devoted to the definition and the properties of this strategy.

\section{Cartesian products of graphs and the fiber-strategy}\label{sect:strat}
In this section, we give the minimal background needed on cartesian products and define an efficient strategy for Gentle, the so-called fiber-strategy, based on the fibers structure of cartesian products of graphs.
For more informations on cartesian product see \cite{sandibook}.
\subsection{Cartesian products of graphs}

Let $H$ and $F$ be two connected simple graphs relatively prime. The vertices of $H\square F$ will be denoted by $(u,v)$, where $u\in V(H)$ and 
$v\in V(F)$. A $H$-fiber is a subgraph of $H\square F$ induced by all the vertices having the same second coordinate. We write $H^v$, where $v\in V(F)$, for 
the $H$-fiber induced by $\{(u,v)|u\in V(H)\}$. Similarly, we define $F^u$, with $u\in V(H)$. 
The $H$-fibers and the $F$-fibers are respectively isomorphic to $H$ and $F$. 
The automorphisms group of $H\square F$ is isomorphic to $Aut(H)\times Aut(F)$.
If $\sigma$ is an automorphism of $H\square F$, it can be seen as a couple $(\psi,\phi)$, where $\psi\in Aut(H)$ and $\phi\in Aut(F)$. In that case, 
$\sigma((u,v))=(\psi(u),\phi(v))$. Another important fact is that $\sigma$ must send a $H$-fiber to another $H$-fiber and the same for the $F$-fibers. 
More precisely, $\sigma(H^v)=H^{\phi(v)}$ and $\sigma (F^u)=F^{\psi (u)}$. To show that a colors preserving automorphism has to be the identity of $Aut(H\square F)$, we will 
mostly proceed as follows. First, we show that an $H$-fiber cannot be sent to another one, which means that $\phi$ is the identity and $\sigma$ fixes the $H$-fibers set wise. 
Using these informations, we prove that $\psi$ is also the identity.

\subsection{Fiber-strategy}

Now, we state some technical results about the fiber-strategy, a strategy that Gentle will follow in mainly all the proofs of our main results.
In a game on $H\square F$, with $H$ non trivial, we say that Gentle follows the \textit{$H$-fiber-strategy}
if we are in one of the following two cases.

\medskip\noindent
{\em Case 0: \begin{itemize}
  \item $|V(H)|$ is even and Rascal starts.
  \item When Rascal plays in a $H$-fiber, Gentle plays in the same $H$-fiber.
 \end{itemize}}

\medskip\noindent
{\em Case 1: 
\begin{itemize}
 \item $|V(H)|$ is odd.
 \item $|V(F)|$ is even and Rascal starts or $|V(F)|$ is odd and Gentle starts.
 \item When Rascal colors the first vertex of a totally uncolored $H$-fiber, Gentle colors the first vertex of another such $H$-fiber.
 \item When Rascal colors a vertex in a $H$-fiber which already has a colored vertex, Gentle colors a vertex in the same $H$-fiber.
\end{itemize}}
The $H$-fiber strategy is always valid.
In Case 0, the parity of each $H$-fiber ensures that Rascal will always be the first to run out of moves in a $H$-fiber. Hence Rascal is always the first to play in each $H$-fiber.
In Case 1, after Gentle's move, there is always an even number of remaining totally uncolored $H$-fibers. Hence Rascal will always be the first to run out of new totally uncolored $H$-fibers to play in.

The following properties are easy and given without proof. There are however fundamental to prove the results in the further sections.

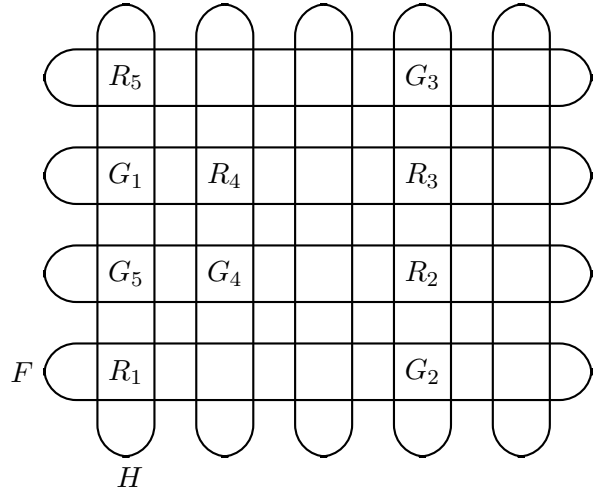
\begin{figure}[!ht]
\begin{center}
\begin{tikzpicture}[scale=1,style=thick]
\def\vr{3pt} 

\def\bunwidth{6}
\def\bunheight{0.75}
\def\spacefactor{1.3}

\foreach \i in {0,1,2,3}
{
	\draw[rounded corners=\bunheight*16] (-\bunwidth/2,\spacefactor * \i-\bunheight/2) rectangle +(\bunwidth * 1.2,\bunheight);
}
\foreach \i in {0,1,2,3,4}
{
	\draw[rounded corners=\bunheight*16] (\spacefactor * \i-\bunwidth/2+0.7,-\bunheight * 1.5) rectangle +(\bunheight,\bunwidth);
        
}
\draw (-\bunwidth/2+\bunheight * 1.5,-\bunheight * 1.5) node[below]{$H$};
\draw (-\bunwidth/2,0) node[left]{$F$};
\draw (-\bunwidth/2+0.7,0) node[right] {$R_1$}; 
\draw (-\bunwidth/2+0.7,\spacefactor * 2-\bunheight/2+0.38) node[right] {$G_1$};
\draw (\spacefactor * 3-\bunwidth/2+0.7,\spacefactor-\bunheight/2+0.38) node[right] {$R_2$}; 
\draw (\spacefactor *3-\bunwidth/2+0.7,-\bunheight/2+0.38) node[right] {$G_2$};
\draw (\spacefactor * 3-\bunwidth/2+0.7,\spacefactor *2-\bunheight/2+0.38) node[right] {$R_3$}; 
\draw (\spacefactor *3-\bunwidth/2+0.7,\spacefactor * 3-\bunheight/2+0.38) node[right] {$G_3$};
\draw (\spacefactor * 1-\bunwidth/2+0.7,\spacefactor *2-\bunheight/2+0.38) node[right] {$R_4$}; 
\draw (\spacefactor *1-\bunwidth/2+0.7,\spacefactor * 1 -\bunheight/2+0.38) node[right] {$G_4$};
\draw ( -\bunwidth/2+0.7,\spacefactor *3-\bunheight/2+0.38) node[right] {$R_5$}; 
\draw ( -\bunwidth/2+0.7,\spacefactor * 1 -\bunheight/2+0.38) node[right] {$G_5$};

\end{tikzpicture}
\end{center}
\caption{How moves alternate in Case 0 of the $H$-fiber strategy (Rascal starts).}
\label{fig:Hstrat}
\end{figure}

\begin{figure}[!ht]
\begin{center}
\begin{tikzpicture}[scale=1,style=thick]
\def\vr{3pt} 

\def\bunwidth{6}
\def\bunheight{0.75}
\def\spacefactor{1.3}

\foreach \i in {0,1,2}
{
	\draw[rounded corners=\bunheight*16] (-\bunwidth/2,\spacefactor * \i-\bunheight/2) rectangle +(\bunwidth * 1.2,\bunheight);
}
\foreach \i in {0,1,2,3,4}
{
	\draw[rounded corners=\bunheight*16] (\spacefactor * \i-\bunwidth/2+0.7,-\bunheight * 1.5) rectangle +(\bunheight, \bunwidth * 0.8);
        
}
\draw (-\bunwidth/2+\bunheight * 1.5,-\bunheight * 1.5) node[below]{$H$};
\draw (-\bunwidth/2,0) node[left]{$F$};
\draw (-\bunwidth/2+0.7,0) node[right] {$G_1$}; 
\draw (-\bunwidth/2+0.7,\spacefactor * 2-\bunheight/2+0.38) node[right] {$G_3$};
\draw (\spacefactor * 3-\bunwidth/2+0.7,\spacefactor * 0-\bunheight/2+0.38) node[right] {$R_1$}; 
\draw (\spacefactor *4-\bunwidth/2+0.7,-\bunheight/2+0.38) node[right] {$G_2$};
\draw (\spacefactor * 3-\bunwidth/2+0.7,\spacefactor *2-\bunheight/2+0.38) node[right] {$R_4$}; 
\draw (\spacefactor *3-\bunwidth/2+0.7,\spacefactor * 1-\bunheight/2+0.38) node[right] {$G_5$};
\draw (\spacefactor * 1-\bunwidth/2+0.7,\spacefactor *2-\bunheight/2+0.38) node[right] {$R_3$}; 
\draw (\spacefactor *2-\bunwidth/2+0.7,\spacefactor * 1 -\bunheight/2+0.38) node[right] {$G_4$};

\draw ( -\bunwidth/2+0.7,\spacefactor * 1 -\bunheight/2+0.38) node[right] {$R_2$};

\end{tikzpicture}
\end{center}
\caption{How moves alternate in Case 1 of the $H$-fiber strategy (Gentle starts).}
\label{fig:Hstrat1}
\end{figure}
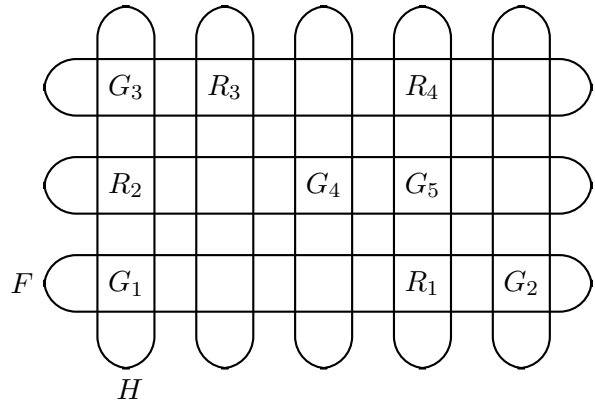

\begin{proposition}\label{prop:strat}
Assume that Gentle plays according to the $H$-fiber-strategy.
\begin{enumerate}
  \item He will color the last vertex of each $H$-fiber.
  \item In Case 0, Rascal will be the first to play in all the $H$-fibers. Then Gentle will play all the second moves in each $H$-fiber, Rascal will play all the third moves and so on.
  \item In Case 1, Gentle will play the first in exactly $\left\lceil \frac {|V(H)|} 2\right\rceil$ different $H$-fibers. Then Rascal will play all the second moves in each $H$-fiber, Gentle will play all the third moves and so on.
\end{enumerate}
\end{proposition}
In Case 0, the moves in a $H$-fiber alternate exactly as in the game played only on $H$, when Rascal starts (see Fig.~\ref{fig:Hstrat}, where $R_i$ and $G_i$ respectively denote the $i^e$ move of Rascal and Gentle). 
This property will be often used by Gentle to play in a $H$-fiber following a winning strategy for the game on $H$.
In Case 1, in a $H$-fiber where Gentle plays first, the moves alternate as in the game played only on $H$, when Gentle starts. 
In a $H$-fiber where Rascal plays first, the only difference is that he is going to play the two first moves in a row (See Fig.~\ref{fig:Hstrat1}, where $R_i$ and $G_i$ have the same meaning as before). 
In that case the lemma beside could be useful. It says that for vertex transitive graphs, if you can win the game playing first, then you can be a real gentleman and let your opponent play this first move.
\begin{lemma}\label{lem:transitive}
 Assume $H$ is vertex transitive. Then either all the first moves are winning for the first player or they are all losing.
\end{lemma}
\proof
 Assume there is $u_0\in V(H)$ such that coloring $u_0$ with $1$ is a winning move for the first player. We have to prove that for any $v\in V(H)$, coloring $v$ with $1$ is also a winning move.
Let $G$ be the game in which the first player has played the winning move $u_0$ and let $G'$ be the game in which his first move has been to color another vertex $v_0$ with $1$. Let $c$ and $c'$ be respectively the coloring built during the game $G$ and $G'$.
Since $H$ is vertex transitive, there exists $\sigma \in Aut(H)$ such that $\sigma(v_0)=u_0$. The winning strategy for the first player in $G'$ is defined by his winning strategy in $G$. When his opponent colors a vertex $w$ in the game $G'$, the first player imagines that his opponent has colored $\sigma(w)$ in $G$ with the same color. There is a vertex $w'$ such that coloring $w'$ is a winning answer for 
the first player in $G$. In the game $G'$, the first player's answer will be to color $\sigma^{-1}(w')$ such that $c'(\sigma^{-1}(w'))=c(w')$.

By assumption, the coloring $c$ is a winning one for the first player. Moreover, for any $v\in V(H)$,  $c(v)=c'(\sigma(v))$.
Hence, an automorphism $\psi$ preserves the coloring $c$ if and only if $\sigma\circ\psi\circ\sigma^{-1}$ preserves the coloring $c'$. This shows that $c'$ is also a winning coloring for the first player. In conclusion, if there is a winning move for the first player, then any first move is a winning move for him.
\qed

\section{Cartesian products of complete graphs.}\label{sec:complete}

Our goal in this section is to prove Theorem~\ref{theo:main-cart} which asserts, under certain conditions on their orders, that for two non trivial connected graphs $H$ and $F$ relatively prime, and of respective order $n$ and $m$, at least one game distinguishing number of $H\square F$ is finite.
Except when both cardinals are equal, it comes directly from the following theorem involving cartesian products of complete graphs. 
Indeed, a distinguishing coloring of $K_{n}\square K_{m}$ is always a distinguishing coloring of $H \square F$. In the first item of Theorem~\ref{theo:main-cart}, when both factors have the same even
cardinal, the corresponding product of complete graphs is not covered by the below result. But, we are in fact going to prove that Gentle's strategy breaks all automorphisms 
of the subgroup of $Aut(K_{n}\square K_{m})$ isomorphic to $Aut(K_{n})\times Aut(K_{m})$. A coloring which distinguishes this subgroup will always be 
a distinguishing coloring of $H\square F$, when the factors are relatively prime.
 
\begin{theorem}\label{theo:complete}
 Let $n$ and $m$ be two distinct natural numbers greater or equal to $2$.
 \begin{enumerate}
  \item If $n\times m$ is even (resp. odd), then $\GDG(K_n\square K_m)=\infty$ (resp.~$\GDR(K_n\square K_m$)).
  \item If $n$ is even, $m\neq n$ and $m\geq n-1$, then $\GDR(K_n\square K_m)\leq m+1$.
  \item If $n$ is odd, $m$ is even and $m\geq 2n-2$, then $\GDR(K_n\square K_m)\leq m+1$.
  \item If $n$ and $m$ are odd and $m\geq 2n-1$, then $\GDG(K_n\square K_m)\leq m+1$.
 \end{enumerate}
\end{theorem}
\proof
The first item is a straightforward application of Proposition~\ref{prop:ordertwo}. For the last items, note that $n\neq m$. Hence, the two factors $K_n$ and $K_m$ are relatively prime.
The vertices of $K_n\square K_m$ are denoted by $(i,j)$, with $i\in\{1,\dots,n\}$ and $j\in\{1,\dots,m\}$. 
The meta-color of a $K_n$-fiber is the list $(c_1,\dots,c_{m+1})$, where $c_l$, with $l\in\{1,\dots,m+1\}$, is the number of vertices in this fiber which are colored with the color $l$ at the end of the game. 
An important remark is that a colors preserving automorphism also preserves the meta-coloring of the $K_n$-fibers.

\begin{figure}[!ht]
\begin{center}
\begin{tikzpicture}[scale=1,style=thick]
\def\vr{3pt} 

\def\bunwidth{7}
\def\bunheight{0.75}
\def\spacefactor{1.6}

\foreach \i in {0,1,2,3}
{
	\draw[rounded corners=\bunheight*16] (-\bunwidth/2,\spacefactor * \i-\bunheight/2) rectangle +(\bunwidth * 1.2,\bunheight);
}
\foreach \i in {0,1,2,3,4}
{
	\draw[rounded corners=\bunheight*16] (\spacefactor * \i-\bunwidth/2+0.7,-\bunheight * 1.5) rectangle +(\bunheight,\bunwidth);
        
}
\draw (-\bunwidth/2+\bunheight * 1.5,-\bunheight * 1.5) node[below]{$K_4$};
\draw (-\bunwidth/2,0) node[left]{$K_5$};
\draw (-\bunwidth/2+0.7,0) coordinate (R1) node[right] {$R_1$}; 
\draw (-\bunwidth/2+0.7,\spacefactor * 2-\bunheight/2+0.38) coordinate (G1) node[right] {$G_1$};
\draw (\spacefactor * 2-\bunwidth/2+0.7,\spacefactor-\bunheight/2+0.38) coordinate (R2) node[right] {$R_2$}; 
\draw (\spacefactor *2-\bunwidth/2+0.7,-\bunheight/2+0.38) coordinate (G2) node[right] {$G_2$};
\draw (\spacefactor * 2-\bunwidth/2+0.7,\spacefactor *2-\bunheight/2+0.38) coordinate (R3)  node[right] {$R_3$}; 
\draw (\spacefactor *2-\bunwidth/2+0.7,\spacefactor * 3-\bunheight/2+0.38) coordinate (G3) node[right]{$G_3$};
\draw (\spacefactor * 4-\bunwidth/2+0.7,\spacefactor *2-\bunheight/2+0.38) node[right] {$R_4$}; 
\draw (\spacefactor *4-\bunwidth/2+0.7,\spacefactor * 1 -\bunheight/2+0.38) node[right] {$G_4$};
\draw ( -\bunwidth/2+0.7,\spacefactor *3-\bunheight/2+0.38)  coordinate(R5) node[right]{$R_5$}; 
\draw ( -\bunwidth/2+0.7,\spacefactor * 1 -\bunheight/2+0.38) coordinate (G5)  node[right]{$G_5$};
\draw ( \spacefactor-\bunwidth/2+0.7,\spacefactor *3-\bunheight/2+0.38)  coordinate(R6) node[right]{$R_6$}; 
\draw ( \spacefactor-\bunwidth/2+0.7,\spacefactor * 0 -\bunheight/2+0.38) coordinate (G6)  node[right]{$G_6$};
\draw ( \spacefactor * 3-\bunwidth/2+0.7,\spacefactor *3-\bunheight/2+0.38) node[right]{$R_7$}; 
\draw ( \spacefactor * 3-\bunwidth/2+0.7,\spacefactor * 1 -\bunheight/2+0.38) node[right]{$G_7$};
\draw ( \spacefactor-\bunwidth/2+0.7,\spacefactor *1-\bunheight/2+0.38)  coordinate(R8) node[right]{$R_8$}; 
\draw ( \spacefactor-\bunwidth/2+0.7,\spacefactor * 2 -\bunheight/2+0.38) coordinate (G8)  node[right]{$G_8$};

\draw(G1) to[bend right] (R1);
\draw(G5) to[bend left] (R5);
\draw(G2) to[bend left] (R2);
\draw(G3) to[bend right] (R3);
\draw(G6) to[bend left] (R6);
\draw(G8) to[bend right] (R8);

\end{tikzpicture}
\end{center}
\caption{Playing according to the matchings in a $K_4$-fiber strategy.}
\label{fig:match}
\end{figure}
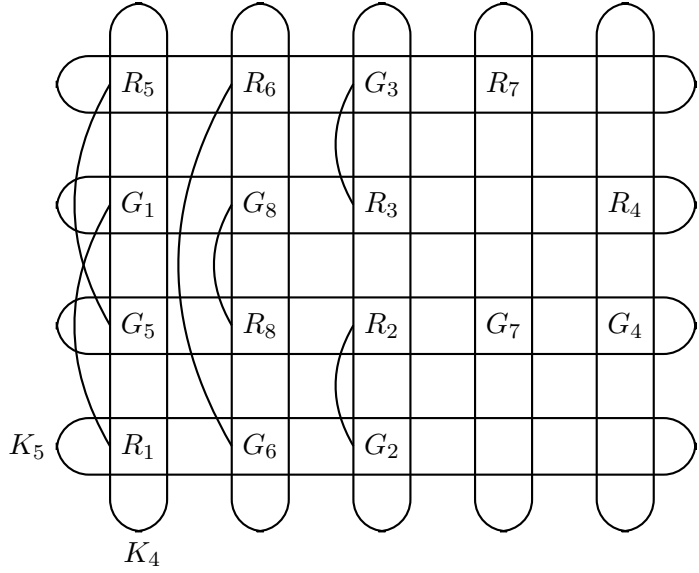

We are going to prove the second statement. We have to give a winning strategy for Gentle playing second with $m+1$ colors.
A proper edge coloring of $K_n$, with $n-1$ colors gives $n-1$ perfect matchings, whose union covers all the edges of $K_n$. We denote these
matchings by $M_1,\dots, M_{n-1}$. Gentle's winning strategy is as follows. First of all, he plays according to the $K_n$-fiber-strategy. 
If Rascal plays in one of the fibers $K_n^j$, with $j\in\{1,\dots,n-1\}$, Gentle plays with respect to the matching $M_j$. It means that
if Rascal colors the vertex $(i,j)$, then Gentle colors the unique vertex $(k,j)$ in $K_n^j$, such that $ik$ is an edge in the matching $M_j$.
Otherwise, he plays as he wants with respect to the $K_n$-fiber-strategy. See Fig.~\ref{fig:match}, where $R_i$ and $G_i$ respectively denote the $i^e$ move of Rascal and Gentle.
Gentle chooses the colors as follows. First, he always plays a color different from the one used by Rascal just before. 
Second, if he has to color the last vertex of a $K_n$-fiber, he chooses the color in a way that the meta-color of this fiber is distinct from all the meta-colors of the already totally colored $K_n$-fibers. He has at most $m-1$ meta-colors to avoid. It is always possible because he can choose among the $m$ colors not used by Rascal just before.

Let us prove now that this strategy yields a distinguishing coloring $c$. Applying the above strategy, Gentle will color the last vertex of each $K_n$-fiber (see Proposition~\ref{prop:strat}). 
Hence, he controls the meta-color of all the $K_n$-fibers. They will all have a distinct meta-color at the end of the game. Therefore, a colors preserving automorphism  $\sigma$ cannot switch these fibers. It means that $\sigma=(\psi,Id_{K_m})$, where $\psi\in Aut(K_n)$.
Assume that $\psi$ is not the identity. There exists $i\in\{1,\dots,n\}$, such that $\psi(i)\neq i$. The edge $\psi(i)i$ belongs to a matching $M_j$, with $j\in\{1,\dots,n-1\}$.
Since Gentle has colored either $(i,j)$ or $(\psi(i),j)$, these two vertices have not the same color. It shows that the automorphism $(\psi,Id_{K_m})$ does not preserve the coloring $c$.
In conclusion $c$ is a distinguishing coloring.

We prove now the third item. The general ideas are similar as above, but since the $K_n$-fibers have odd order, a matching does not cover all the vertices of $K_n$.
Hence, if Rascal is the first to play in a $K_n$-fiber, Gentle cannot immediately play with respect to a matching. 
Since $n$ is odd, we need $n$ matchings to have that their union covers all the edges of $K_n$. Each matching does not cover exactly one vertex and this uncovered vertex is distinct for each of them. We denote by $M_j$, with $j\in\{1,\dots,n\}$, the matching which does not cover the vertex $j$. 
Without lost of generality, we assume that Rascal's first move is to color the vertex $(1,1)$ in $K_n^1$. Gentle will again follow a $K_n$-fiber strategy. 
Hence, he will be the first to color a vertex in exactly $\frac m 2$ different $K_n$-fibers, say $K_n^2,\dots,K_n^{\frac m 2+1}$. 
When Gentle colors a vertex of $K_n^j$, with $j\in\{1,\dots,\frac{m}{2}+1\}$, if it is the first vertex of $K_n^j$ to be colored, he chooses the vertex $(j,j)$.
Otherwise, he plays with respect to the matching $M_j$. When the fiber-strategy leads him to play in other $K_n$-fibers, he plays wherever he wants with respect to the $K_n$-fiber strategy. 
For the choice of the colors, he plays as in the previous strategy. The proof that the coloring built during the game is distinguishing is exactly the same. 
Just note that by hypothesis, $\frac m 2+1\geq n$. Hence, Gentle has enough $K_n$-fibers to use the $n$ matchings needed to cover all the edges of $K_n$. 

The proof of the last item is the same as for the previous one. Because Gentle starts, he will be the first to play in $\frac{m+1}{2}$ different $K_n$-fibers, which is by hypothesis greater or equal to $n$.
\qed

For $K_2\square K_m$, we can compute the exact value of $\GDR$. In \cite{schmidt}, it is shown that $\GDR(K_2\square K_m)=3$, for $m\in \{2,3,4\}$. 
For $m\geq 5$, we are going to prove that we need exactly $m$ colors. Then, the bound obtained above is not far to be tight. 

Let $c$ be a coloring of $K_2\square K_m$, with $m\geq 5$. We say that two distinct $K_2$-fibers, $K_2^i$ and $K_2^j$ are colored the same if $c(K_2^i)=c(K_2^j)$.
If we have also $c((i,1))=c((j,1))$, we say that the two fibers are strictly colored the same. 
\begin{proposition}\label{prop:Kn<=}
If $m\geq 5 $, then $\GDR(K_2\square K_m)= m$.
\end{proposition}
\proof
First, we show that with $m$ distinct colors Gentle has a winning strategy. Recall that Rascal starts. 
When Rascal plays in a $K_2$-fiber, Gentle answers by coloring the second vertex of this $K_2$-fiber. That means he plays according to a $K_2$-fiber strategy. 
He colors in a way that the new colored $K_2$-fiber is not colored the same as another $K_2$-fiber already colored before.
This is always possible, because there are at most $m-1$ different $K_2$-fibers colored before and Gentle can use $m$ colors. Moreover, he can ensure that at least one $K_2$-fiber is not monochromatic.
Let us prove now that this strategy yields a distinguishing coloring. Assume $\sigma$ is a colors preserving automorphism. Then $\sigma(K_2^i)=K_2^i$, for all $i\in\{1,\dots,m\}$. 
But, there is at least one bi-chromatic $K_2$-fiber. Hence, $\sigma$ must also fix this $K_2$-fiber point wise. Therefore, $\sigma$ has to fix all the $K_2$-fibers point wise. In conclusion, $\sigma$ is the identity. 

It remains to prove that Rascal has a winning strategy, if strictly less than $m$ colors are allowed during the game.
Remark that, if two distinct $K_2$-fibers, $K_2^i$ and $K_2^j$ are strictly colored the same at any moment in the game, then Rascal wins. 
Indeed, there is an automorphism $\sigma$ such that $\sigma((i,1))=(j,1)$, $\sigma((i,2))=(j,2)$ and $\sigma$ fixes all the other vertices.

\noindent
Rascal starts by coloring $(1,1)$ with $1$. There are two cases.\\
\textbf{Case 1:} Gentle colors the vertex $(1,2)$. 

\noindent
Rascal answers by coloring  $(2,1)$ with $1$.
If Gentle colors a vertex different than $(2,2)$, Rascal wins at his next turn by coloring $(2,2)$ with the same color as $(1,2)$. So, we can suppose that Gentle colors $(2,2)$. Turn by turn, this shows that Rascal can color all the vertices of the form $(i,1)$ with the color $1$, and that Gentle is forced
to color only the vertices of the form $(i,2)$. Since Gentle has strictly less than $m$ colors at his disposal, there are two vertices $(i_0,2)$, $(j_0,2)$, which receive the same color. 
Hence, the two $K_2$-fibers, $K_2^{i_0}$ and $K_2^{j_0}$ will be strictly colored the same and Rascal will win.

\noindent
\textbf{Case 2:} Gentle first move is to color the vertex $(2,x)$, with $x\in\{1,2\}$. 

\noindent
Rascal answers by coloring with $1$ the vertex
$(3,1)$. Now, if Gentle plays in $K_2^1$ or $K_2^3$, Rascal wins because he can play such that $K_2^1$ and $K_2^3$ are strictly colored the same. 
Suppose that Gentle plays a vertex which is not in $K_2^1$ or $K_2^3$. Since $m\geq 5$, at least one vertex in the fiber $K_m^1$ is still uncolored, say $(5,1)$. Rascal replies by coloring this vertex with $1$. The vertices $(1,2)$, $(3,2)$ and $(5,2)$ are still uncolored. Rascal can ensure that at least two of the three $K_2$-fibers, $K_2^1$, $K_2^3$ and $K_2^5$ are strictly colored the same. 
Indeed, if Gentle is the first to color one of these three uncolored vertices, Rascal copies this color in one of the two remaining vertices. 
Otherwise, he will be able to decide the coloring of two of them. In conclusion, Rascal will also win in this second case.
\qed

Of course, Theorem~\ref{theo:complete} does not cover all possibilities. We did not manage to prove that in the remaining cases the invariants are finite. 
But we know that the $K_n$-fiber strategy used above by Gentle will fail in these cases. More precisely, we have the following proposition.
\begin{proposition} Let $n$ and $m$ be two distinct natural numbers greater or equal to $2$. Whatever the number of colors allowed, if Gentle follows a $K_n$-fiber strategy on $K_n\square K_m$, he looses in both following cases:
\begin{itemize}
 \item Rascal starts, $n$ is odd, $m$ is even and $m<2n-2$,
 \item Gentle starts, $n$ and $m$ are odd and $m<2n-1$.
\end{itemize}
\end{proposition}
\proof
 We prove the first statement. The second can be proved in exactly the same way. Rascal winning strategy is to create two $K_m$-fibers, say $K_m^1$ and $K_m^n$, which are strictly colored the same. 
More precisely, if $u\in K_m^1$ and $v\in K_m^n$ are in the same $K_n$-fiber then they share the same color. 
In that case, the automorphism which only permutes $K_m^1$ and $K_m^n$ is a colors preserving automorphism. 

Rascal proceeds as follows. He plays his $\frac m 2$ first moves in the same $K_m$-fiber, say $K_m^1$. Since Gentle plays according to a $K_n$-fiber strategy, at the end of the $(\frac m 2)^\text{th}$ turn of Gentle each $K_n$-fiber has exactly one colored vertex. 
These $m$ first moves are called the first phase of the game. Let $k$ be the number of uncolored vertices in $K_m^1$ at the end of this first phase. We have: $0\leq k\leq \frac m 2$ ($k$ could be equal to $0$, if Gentle has only played in $K_m^1$ during the first phase). 
The forthcoming $k$ moves of Rascal and $k$ moves of Gentle will be called the second phase of the game. In this phase, when Rascal plays in a $K_n$-fiber, Gentle has to answer by a move in this same $K_n$-fiber. Hence, Rascal can play all the $k$ remaining uncolored vertices of $K_m^1$. 
Let $u$ be such a vertex. There is a unique vertex $v$ in the same $K_n$-fiber than $u$, which is already colored (this vertex has been colored by Gentle during the first phase). Rascal copies the color of $v$ to color $u$. 
During this second part of the game, Gentle has played in at most $k$ distinct $K_m$-fibers. 
Since $m<2n-2$, then $k<n-1$. Hence, there exists a $K_m$-fiber, say $K_m^n$, in which Gentle has not played during this second phase. In $K_m^n$, there is at most one colored vertex,
say $w$. In that case, $w$ has been colored by Gentle during the first phase of the game. 
This vertex $w$ shares the same color as the vertex of $K_m^1$, which is in the same $K_n$-fiber (this vertex has been colored by Rascal during the second phase). 
Therefore, Rascal can now color all the uncolored vertices of $K_m^n$, such that $K_m^1$ and $K_m^n$ are strictly colored the same. 
\qed

\section{Cartesian products of involutive graphs.}\label{sect:inv}

In this section, we study the game distinguishing numbers of cartesian products of involutive graphs and prove Theorem~\ref{theo:main-inv}.
The class of involutive graphs has been introduced in \cite{schmidt}. It contains graphs like even cycles, hypercubes or more generally diametrical graphs and even graphs (see \cite{Gobel}). 
Let us recall the definition. An {\em involutive graph} $H$ is a graph together with an involution, $bar:V(H)\rightarrow V(H)$, which commutes with all automorphisms
and has no fixed point. In other words:\begin{itemize}
                 \item $\forall u\in V(H)$, $ \overline{\overline u}=u$ and $\overline u\neq u$,\\
                 \item $\forall \sigma\in Aut(H)\; \forall u\in V(H)$, $\sigma(\overline u)=\overline{\sigma(u)}$.\\
                 \end{itemize}
The set $\{u,\overline u\}$ will be called a {\em block}.
An important remark is that an automorphism of an involutive graph has to map a block to a block. In other words, there is a natural action of the automorphism group on the set of blocks.
We introduce the following concepts, which are going to play a similar role as the meta-colors used in Section~\ref{sec:complete}.

If $H$ is an involutive graph and $c$ is a vertex-coloring with $d$ colors, then the {\em type} $t$ of a block $\{u,\overline u\}$ is defined by:

\noindent
\begin{center}
\begin{math}
t(\{u,\overline u\})=\begin{cases}
                    c(u)-c(\overline u) \bmod d &\text{if } c(u)-c(\overline u)\bmod d\in\{0,\dots,\lfloor\frac d 2\rfloor\}\\
                    c(\overline u)-c(u)\bmod d &\text{otherwise.}
                   \end{cases}
\end{math}
\end{center}

\noindent
The {\em block-list} of $H$, $L_c(H)$ is the list $(n_0,\dots,n_{\lfloor\frac d 2\rfloor})$ of length $\lfloor\frac d 2\rfloor+1$, where $n_i$, with $i\in\{0,...\lfloor\frac d 2\rfloor\}$, is the number
of blocks of type $i$, according to the coloring $c$. Note that, if $\sigma$ is an automorphism of $H$, then $t(\sigma(\{u,\bar u\}))=t(\{u,\bar u\})$ and $L_c(H)=L_{c}(\sigma(H))$. 

Assume now that $H$ is a connected involutive graph and $F$ is a connected graph relatively prime to $H$. The following proposition asserts that if the
classical distinguishing number of $F$ is not too big, then $\GDR(H\square F)$ is bounded above by $\GDR(H)$. Theorem~\ref{theo:main-inv} will be a straightforward application
of this result.  
\begin{theorem}\label{theo:inv}
Let $H$ be a connected involutive graph. Assume that Gentle has a winning strategy playing second on $H$, with $d\geq\GDR(H)$ colors. 
Moreover, this strategy yields colorings, whose block-list is in a fixed set $\mathcal L$ .

\medskip\noindent
If $F$ is a connected graph relatively prime to $H$, with $D(F)\leq  \left(\begin{array}{c}
                \frac{|V(H)|}2+\lfloor\frac d 2\rfloor\\\
                \lfloor\frac d 2\rfloor\
              \end{array}\right)-|\mathcal{L}|+1$,
then $\GDR(H\square F)\leq d$.
\end{theorem}
\proof 
We have to give a Gentle winning strategy with $d$ colors, assuming Rascal starts. The coloring obtained at the end of the game will be denoted by $c$.

First of all Gentle will follow a $H$-fiber strategy. Note that $|V(H)|$ is even. Hence, we are in Case $0$ of this strategy. Let $(u_1,v_1)$ be the first vertex of $H\square F$ colored by Rascal. 
Gentle imagines a distinguishing coloring $c'$ of $F$, with $D(F)$ colors. When Gentle has to play in the fiber $H^{v_1}$, he chooses the vertex and the color according to a winning strategy in $H$. 
This is possible, because 
Gentle's moves and Rascal's moves in $H^{v_1}$ alternate like the moves in a game played only on $H$, when Rascal starts (see Proposition~\ref{prop:strat}). Moreover, $d\geq \GDR(H)$ by hypothesis.
In the other $H$-fibers, when Rascal plays the vertex $(u,v)$, Gentle answers  by coloring the vertex $(\overline u,v)$. In this way, Gentle will be able to control 
the block-list of these fibers. More precisely, he chooses the colors such that:
$$(\dag)~~\forall v,w\in V(F),~L_c(H^v)=L_c(H^w) \text{ only if } c'(v)=c'(w).$$
This is possible if there exists at least $(D(F)-1)+|\mathcal {L}|$ distinct possible block-lists. Indeed, Gentle cannot control in advance the block-list of the fiber $H^{v_1}$. 
By hypothesis, we only know that this block-list will belong to $\mathcal L$. Hence, $|\mathcal L|$ kinds of block-list are used to stand for the imaginary color $c'(v_1)$. 
Finally, with $(D(F)-1)+|\mathcal {L}|$ possible block-lists, Gentle has enough possibilities to associate distinct block-lists to distinct colors of the coloring $c'$. 
The number of block-lists is the number of weak compositions of $\frac {|V(H)|}2$ (the number of blocks) into $\lfloor\frac d 2\rfloor+1$ 
natural numbers (the number of block types). So, there are $\left(\begin{array}{c}
                \frac{|V(H)|}2+\lfloor\frac d 2\rfloor\\\
                \lfloor\frac d 2\rfloor\
              \end{array}\right)$ kinds of block-lists, which is by hypothesis greater or equal to $(D(F)-1)+|\mathcal L|$. 

Now, we prove that the coloring obtained with this strategy is distinguishing. Assume $\sigma$ is a colors preserving automorphism. 
We have $\sigma=(\psi,\phi)$, where $\psi\in Aut(H)$ and $\phi\in Aut(F)$. This automorphism maps blocks in $H^{v}$ to blocks in $H^{\phi(v)}$, for any $v\in V(F)$. Hence, the automorphism $\phi$ preserves the block-lists of the $H$-fibers: $L_c(H^v)=L_c(H^{\phi(v)})$, for all $v\in V(F)$. By condition $(\dag)$, this automorphism preserves also the distinguishing coloring $c'$ of $F$. 
Hence, $\phi$ is the identity of $Aut(F)$. This implies that $\sigma(H^{v_1})=H^{v_1}$. 
But the coloring of this $H$-fiber is obtained by following a winning strategy for Gentle in the game on $H$. Therefore, $\psi$ is the identity of $H$. In conclusion, $\sigma$ is trivial and the coloring $c$ is distinguishing.  
\qed

Theorem~\ref{theo:main-inv} is a straightforward application of the above result for two reasons. First, we know that  for an involutive graph $H$, $\GDR(H)\leq D^2(H)+D(H)-2$ 
(see Theorem~\ref{theo:finite-inv}). 
Moreover, with  $D^2(H)+D(H)-2$ colors, Gentle has a winning strategy such that he knows exactly the block-list he will get at the end of the game (see the proof of Theorem~1.6 in \cite{schmidt}). It means, with the notation of the above theorem, that $\mathcal L$ is just a singleton. 

\section{Cartesian products of cycles}\label{sec:tore}

In this final section, we give a proof of Theorem~\ref{theo:main-tore}. Note that proving the statement about the infinity of the invariants is a straightforward application of Proposition~\ref{prop:ordertwo}. 
For the cycle $C_n$ of order $n\geq 3$, we set $V(C_n)=\{1,\dots,n\}$ and $E(C_n)=\{ij\;|\;|i-j|=1\bmod n,\;i,j\in V(C_n)\}$.
We begin with toroidal grids of even order. Since even cycles are involutive graphs the first proposition is a direct corollary of Theorem~\ref{theo:inv}.

\begin{proposition}\label{prop:even} Let $n$ and $m$ be two distinct natural numbers greater or equal to $3$.
If $n$ is even and $n\geq 8$, then $\GDR(C_n\square C_m)=2$.
\end{proposition}

\proof
 In \cite[Proposition~4.1]{schmidt}, the winning strategy used by Gentle with two colors leads to exactly three bi-chromatic blocks, when $n\geq 12$.
To show that $\GDR(C_8)=\GDR(C_{10})=2$, they used an exhaustive computer check. This computing also gives that there is a Gentle's winning strategy which leads to one or three bi-chromatic blocks if $n=8$, and to one or four bi-chromatic blocks if $n=10$. 
Therefore, with the same notations as in Theorem~\ref{theo:inv}, we have that $|\mathcal L|\leq 2$. For all $m\geq 3$, $D(C_m)\leq 3$. Thus, we have
$D(C_m)\leq \left(\begin{array}{c}
                \frac{n+2}2\\
                1
              \end{array}\right)-2+1$, and we can directly applied Theorem~\ref{theo:inv} to get that $\GDR(C_n\square C_m)=2$.
\qed
\begin{proposition}\label{prop:C4C6} Let $n$ be in $\{4,6\}$ and $F$ a connected graph relatively prime to $C_n$, with at least three vertices.
If $D(F)\leq 3$, then $\GDR(C_n\square F)=2$.
\end{proposition}
\proof
Let $n\in\{4,6\}$. We denote by $c$ the coloring built during the game. We have to give a winning strategy for Gentle with $2$ colors. Gentle plays according to a $C_n$-fiber strategy and uses the block-lists as meta-colors. Here, the problem is that $\GDR(C_n)=3$.
Gentle fancies a distinguishing coloring $c'$ of $F$, where the three colors are really used. 
As in Theorem \ref{theo:inv}, Gentle can control the block-list of the $C_n$-fibers such that:
$$\forall v\in V(F),~L_c(C_n^v)=\begin{cases}
                                       (n,0)&\text{ if } c'(v)=1,\\
                                       (n-1,1)&\text{ if } c'(v)=2,\\
                                       (n-2,2)&\text{ if } c'(v)=3.\\
                                      \end{cases}$$
Moreover, for $v\in V(F)$, if $L_c(C_n^v)$ must be equal to $(n-1,1)$ or $(n-2,2)$, he plays such that the block $\{(1,v),(n/2,v\}$ is of type $1$. 

Now, we prove that the coloring $c$ is distinguishing. Assume $\sigma$ is a colors preserving automorphism. We have $\sigma=(\psi,\phi)$, where $\psi\in Aut(C_n)$ and $\phi\in Aut(F)$. 
For all $v\in V(F)$, $L_c(C_n^v)=L_c(C_n^{\phi(v)})$. Hence, for all $v\in V(F)$, $c'(v)=c'(\phi(v))$. Since $c'$ is a distinguishing coloring of $F$, we get that $\phi$ is trivial. 
Hence, $\sigma$ fixes the $C_n$-fibers set wise. Since there is at least one $C_n$-fiber with block-list $(n-1,1)$, $\psi$ must be the symmetry $\Delta$ of axes $(1,\frac n 2)$ or the identity. 
But $\Delta$ does not preserve the coloring in the $C_n$-fibers, whose block-list is $(n-2,2)$.
Indeed, in such a fiber, one of the block of type $1$ is stable under $\Delta$. The other block of type $1$ is sent by $\Delta$ to a block of type $0$ or switched to itself. In both cases, it breaks the coloring. 
In conclusion $\psi$ is the identity and so is $\sigma$. \qed

This result directly implies the following corollary.
\begin{corollary}\label{cor:C4C6} Let $m$ be an integer greater or equal to $3$.
\begin{enumerate}
 \item If $m\neq 6$, then $\GDR(C_6\square C_m)=2$.
 \item If $m\neq 4$, then $\GDR(C_4\square C_m)=2$.
\end{enumerate}
\end{corollary}
The following proposition gives a general upper bound, when one factor has distinguishing number less or equal to $2$. It has as corollary, the case where both factors have odd cardinal and a least one is not prime.
\begin{proposition}\label{prop:cart2}
Let $H$ and $F$ be two connected graphs relatively prime. Assume $H$ is vertex transitive, $D(H)\geq 2$ and $D(F)\leq2$.
	    \begin{enumerate}
	    \item If $|V(H)|$ and $|V(F)|$ are odd, then $\GDG(H\square F)\leq\GDG(H)$.
	     \item If $|V(H)|$ is odd and $|V(F)|$ is even, then $\GDR(H\square F)\leq\GDG(H)$.
	     \item If $|V(H)|$ is even, then $\GDR(H\square F)\leq\GDR(H)$
	    \end{enumerate}
\end{proposition}
\proof We prove the first statement. Let $c$ be the coloring built during the game. For each $H$-fiber $H^v$, with $v\in V(F)$, we define: 
$$p(H^v)=\begin{cases}
           1&\text{ if } |\{u\in H^v|c((u,v))=1\}| \text{ is odd}\\
           2&\text{ otherwise.}
          \end{cases}$$
We have to define a Gentle's winning strategy with $\GDG(H)$ colors. Gentle is going to play according to a $H$-fiber strategy. Note that we are in Case 1 of this strategy.  In the $H$-fibers, where Gentle is the first to play, the moves alternate exactly as in a game played only on $H$, with Gentle playing first (see Proposition \ref{prop:strat}). 
In the other $H$-fibers, it is also the case, except for the first move which is played by Rascal. 
In other words, Rascal will play the two first moves in a row in these $H$-fibers. 
Since $H$ is vertex transitive, we assume, by Lemma~\ref{lem:transitive}, that Gentle has played first also in these $H$-fibers. 
Therefore, Gentle can play following a winning strategy for $H$ in each $H$-fiber. He plays like this as long as one $H$-fiber is totally colored,
say $H^{v_0}$. Now, he fancies a distinguishing coloring $c'$ of $F$ such that $c'(v_0)=p(H^{v_0})$. For the later moves, he plays such that: 
$$(\ddagger)\;\; \forall v\in V(F),\; c'(v)=p(H^v).$$
Since he follows a $H$-fiber strategy, we recall that he is going to play the last move in each $H$-fiber. Hence, he is able to decide the parity of the number of vertices colored with $1$ in each of them.

Let us prove now that the coloring $c$ is distinguishing. Let $\sigma$ be a colors preserving automorphism. For all $v\in V(F)$, we have $p(\sigma(H^v))=p(H^v)$. Since $c'$ is a distinguishing coloring of $F$, it implies, by $(\ddagger)$, 
that $\sigma$ fixes the $H$-fibers set wise. Therefore, $\sigma(H^{v_0})=H^{v_0}$. But the coloring on this $H$-fiber is obtained by following a winning strategy on $H$. Then, $H^{v_0}$ must be fixed point wise by $\sigma$. In conclusion, $\sigma$ is the identity.

For the two remaining statements, the proof is almost the same. The only difference is that for the third item, we are in Case 0 of the $H$-fiber strategy. \qed

\begin{corollary}\label{cor:odd}
 Let $n$ and $m$ be two odd distinct natural numbers greater or equal to $3$. If $n$ is not prime and $m\geq7$, then $\GDG(C_n\square C_m)=2$.
\end{corollary}
\proof
 Under the hypothesis of the corollary, we have $\GDG(C_n)=2$ and $D(C_m)=2$. Thus, this is a straightforward application of Proposition~\ref{prop:cart2}. 
\qed

With the previous results, we are able to compute the distinguishing numbers of the toroidal grid $C_n\square C_m$, except for the following cases:
\begin{itemize}
\item $C_3\square C_m$, with $m\neq 3$ and $m$ odd,
\item $C_5\square C_m$, with $m\neq 5$ and $m$ odd,
\item $C_n\square C_m$, with $n\neq m$, $n$ and $m$ odd and prime.
\end{itemize}
To settle this remaining cases, we state the following proposition.
\begin{proposition}\label{prop:prime}
 Let $n$ and $m$ be two distinct odd numbers greater or equal to $3$. If $n$ is prime and $m\geq 7$, then $\GDG(C_n\square C_m)=2$
\end{proposition}
\proof
 Let $c$ be the coloring built during the game. For each $C_n$-fiber $C_n^j$, with $j\in\{1,...m\}$, we define: 
$$p(C_n^j)=\begin{cases}
           1&\text{ if } |\{i\in C_n^j|c((i,j))=1\}| \text{ is odd}\\
           2&\text{ otherwise.}
          \end{cases}$$
Let $M_1$, $M_2$ and $M_3$ be three maximum matchings of $C_n$, whose union covers $E(C_n)$. Let $v_1$, $v_2$ and $v_3$ be the only vertices of $C_n$, which are respectively not covered by $M_1$, $M_2$ and $M_3$. 
Let $c'$ be a distinguishing coloring of $C_m$, with $2$ colors. Such a coloring exists because $m>5$.

We have to outline a Gentle's winning strategy with $2$ colors. He is going to follow a $C_n$-fiber strategy. 
Since $m\geq7$,  there are at least three distinct $C_n$-fibers,  $C_n^{j_1},C_n^{j_2},C_n^{j_3}$, with $j_1,j_2,j_3\in\{1,\dots,m\}$,  where Gentle is the first to play. 
The first vertex that Gentle is going to color in $C_n^{j_k}$, with $k\in\{1,2,3\}$ is $(v_k,j_k)$. He colors it such that $c((v_k,j_k))=c'(j_k)$.
In the other $C_n$-fibers, where he is the first to play, the first vertex he chooses and the color he uses do not matter. For the later moves in $C_n^{j_k}$, he will choose the vertices with respect to the matching $M_k$. 
Moreover, he uses the other color than the one used by Rascal just before. In this way, the parity of the number of vertices in $C_n^{j_k}$ colored with $1$ only depends on $c((v_k,j_k))$. 
Hence $p(C_n^{j_k})=c((v_k,j_k))=c'(j_k)$, for $k\in\{1,2,3\}$.
For the moves in $C_n^j$, with $j\not\in\{j_1,j_2,j_3\}$, Gentle plays whatever he wants, except when he colors the last vertex of the $C_n$-fiber. In that case, he chooses 
the color such that $p(C_n^j)=c'(j)$.

We prove now that $c$ is distinguishing. Let $\sigma$ be a colors preserving automorphism. For all the $C_n$-fibers, we have $p(C_n^j)=c'(j)$. Since $c'$ is a distinguishing coloring of $C_m$, $\sigma$ fixes the $C_n$-fibers set wise. 
Thus, we have $\sigma=(\psi,Id)$, with $\psi\in Aut(C_n)$. Since $n$ is prime, any non trivial rotation acts transitively on a $C_n$-fiber.
As at least one such fiber is not monochromatic, $\psi$ could not be a non trivial rotation. In the other hand, if $\psi$ is an axial symmetry, since
 $n$ is odd, there is an edge $e\in E(C_n)$ such that $\psi(e)=e$. This edge belongs to one of the three matchings, say $M_1$. Gentle has played such that in $C_n^{j_1}$, the edge corresponding to $e$ is not monochromatic. Therefore $\psi$ cannot preserve the coloring. In conclusion, $\psi$ must be the identity and so is $\sigma$. 
\qed

We are now ready to prove Theorem~\ref{theo:main-tore}.

\noindent
\textbf{Proof of Theorem~\ref{theo:main-tore}.} Let $n_1,\dots,n_k$, with $k\geq2$, be $k$ distinct numbers greater or equal to $3$. 
If $\displaystyle\prod_{i=1}^k n_i$ is even (resp. odd), we have to prove that $\GDR(C_{n_1}\square\cdots\square C_{n_k})=2$ (resp. $\GDG(C_{n_1}\square\cdots\square C_{n_k})=2$).
If $k=2$, this is a consequence of Propositions~\ref{prop:even} and\ref{prop:prime} and Corollaries~\ref{cor:C4C6} and \ref{cor:odd}, except for $C_3\square C_5$. 
An exhaustive computer check prove that in that case two colors are also enough. 
For $k\geq 3$, we proceed by induction. If we are not dealing with $C_3\square C_4\square C_5$, we can assume that $n_k\geq 6$. 
Hence, $D(C_{n_k})=2$ and by induction $\GDR(C_{n_1}\square\cdots\square C_{n_{k-1}})=2$ or $\GDG(C_{n_1}\square\cdots\square C_{n_{k-1}})=2$, depending on the parity. 
Finally, we apply Proposition~\ref{prop:cart2}, with $H=C_{n_1}\square\cdots\square C_{n_{k-1}}$ and $F=C_{n_k}$, to get the expected results. For $C_3\square C_4\square C_5$, we apply Proposition~\ref{prop:C4C6} to show that $\GDR(C_3\square C_4\square C_5)=2$.  \qed

\begin{corollary}
 Let $n_1,\dots,n_k$ be $k$ distinct natural numbers greater or equal to $2$. 
\begin{enumerate}
\item If $\displaystyle\prod_{i=1}^k n_i$ is even, then $\GDG(P_{n_1}\square\cdots\square P_{n_k})=\infty$ and $\GDR(P_{n_1}\square\cdots\square P_{n_k})=2$.
 \item If $\displaystyle\prod_{i=1}^k n_i$ is odd, then $\GDR(P_{n_1}\square\cdots\square P_{n_k})=\infty$ and $\GDG(P_{n_1}\square\cdots\square P_{n_k})=2$.
 \end{enumerate}
\end{corollary}
\proof If $k=1$, we easily have $\GDR(P_n)=2$, when $n$ is even and $\GDG(P_n)=2$, when $n$ is odd (see \cite{schmidt}).
If $k\geq 2$ and $n_i\geq 3$, for all $i\in\{1,\dots,k\}$, then it is a straightforward consequence of Theorem~\ref{theo:main-tore}. 
Indeed, in this case, a distinguishing coloring of $C_{n_1}\square\cdots\square C_{n_k}$ is also a distinguishing coloring of $P_{n_1}\square\cdots\square P_{n_k}$. If one factor, say $P_{n_1}$ is isomorphic to 
$P_2$, then we can apply Proposition~\ref{prop:cart2}, with $H=P_{n_1}$ and $F=P_{n_2}\square\cdots\square P_{n_k}$. We have actually that $D(P_{n_2}\square\cdots\square P_{n_k})=2$ (the only cartesian products of paths for which it is not true are $P_2\square P_2$ and $P_2\square P_2\square P_2$ , with $D(P_2\square P_2)=D(P_2\square P_2\square P_2)=3$).
\qed

\section*{Acknowledgements}

The research was in part financed by the ANR-14-CE25-0006 project of the French National Research Agency.


\begin{thebibliography}{}
\bibitem{albertson} M.O.~Albertson and K.L.~Collins, Symmetry breaking in graphs, Electron.~J.~Comb. 3 (1996), \#R18.

\bibitem{bogstad} B.~Bogstad  and L.~Cowen, The distinguishing number of hypercubes, Discrete Math. 383 (2004) 29--35.


\bibitem{faigle} U.~Faigle, U.~Kern, H.~Kierstead and W.T.~Trotter, On the game chromatic number of some classes of graphs, Ars Combinatoria 35 (1993) 143--150.

\bibitem{fisher} M.J.~Fisher and G.~Isaak, Distinguishing colorings of cartesian products of complete graphs, Discrete Math. 308 (2008) 2240--2246

\bibitem{Gobel} F.~G\"obel and H.J.~Veldman, Even graphs, J.~Graph Theor. 10 (1986) 225--239.

\bibitem{klav.cart.clik} W.~Imrich, J.~Jerebic and S.~Klav\v{z}ar, The distinguishing number of cartesian products of complete graphs, Eur.~J.~Comb. 45 (2009) 175--188.

\bibitem{Imrich_cartes_power} W.~Imrich and S.~Klav\v{z}ar, Distinguishing cartesian powers of graphs, J.~Graph Theor. 53 (2006) 250--260.

\bibitem{sandibook} W.~Imrich and S.~Klav\v{z}ar, Product graphs: struture and recognition, Wiley-Interscience Series in Discrete Mathematics and Optimization (2000).

\bibitem{klwozh-06} S.~Klav\v zar, T.L.~Wong and X.~Zhu, Distinguishing labelings of group action on vector spaces and graphs, J.~Algebra 15 Issue 2 (2006) 626--641.

\bibitem{klav_power} S.~Klav\v{z}ar and X.~Zhu, Cartesian powers of graphs can be distinguished by two labels, Eur.~J.~Comb. 28 (2007) 303--310.


\bibitem{schmidt} S.~Gravier, K.~Meslem, S.~Schmidt and S.~Slimani, A New game invariant of graphs: the game distinguishing number, Preprint arXiv:1410.3359v4 [math.CO] (2015).

\bibitem{domgame}  
B.~Bre\v{s}ar, S.~Kla\v{z}ar and D.F.~Rall, Domination game and an imagination strategy, SIAM J.~Discrete Math. 24 (2010) 979--991. 



\end{thebibliography}
\end{document}